\documentclass{amsart}
\usepackage{amssymb}
\usepackage{amsfonts}

\newtheorem{theorem}{Theorem}
\theoremstyle{plain}

\newtheorem{corollary}{Corollary}

\newtheorem{definition}{Definition}

\newtheorem{lemma}{Lemma}

\newtheorem{remark}{Remark}

\numberwithin{equation}{section}
\input{tcilatex}

\begin{document}
\title[Weighted Morrey estimates]{Some inequalities for the multilinear
singular integrals with Lipschitz functions on weighted Morrey spaces}
\author{FER\.{I}T G\"{U}RB\"{U}Z}
\address{ Hakkari University, Faculty of Education, Department of
Mathematics Education, Hakkari 30000, Turkey}
\email{feritgurbuz@hakkari.edu.tr}
\urladdr{}
\thanks{}
\curraddr{ }
\urladdr{}
\thanks{}
\date{}
\dedicatory{}
\thanks{}

\begin{abstract}
The aim of this paper is to prove the boundedness of the oscillation and
variation operators for the multilinear singular integrals with Lipschitz
functions on weighted Morrey spaces.
\end{abstract}

\maketitle

\section{Introduction}

\bigskip We first say that there exists a continuous function $K\left(
x,y\right) $ defined on $\Omega =\left\{ \left( x,y\right) \in {\mathbb{%
R\times R}}:x\neq y\right\} $ and $C>0$ if $K$ admits the following
representation 
\begin{equation}
\left\vert K\left( x,y\right) \right\vert \leq \frac{C}{\left\vert
x-y\right\vert },\qquad \forall \left( x,y\right) \in \Omega  \label{1}
\end{equation}%
and for all $x$, $x_{0}$, $y\in {\mathbb{R}}$ with $\left\vert
x-y\right\vert >2\left\vert x-x_{0}\right\vert $ 
\begin{eqnarray}
&&\left\vert K\left( x,y\right) -K\left( x_{0},y\right) \right\vert
+\left\vert K\left( y,x\right) -K\left( y,x_{0}\right) \right\vert  \notag \\
&\leq &\frac{C}{\left\vert x-y\right\vert }\left( \frac{\left\vert
x-x_{0}\right\vert }{\left\vert x-y\right\vert }\right) ^{\beta },
\label{1*}
\end{eqnarray}%
where $1>\beta >0$. Then $K$ is said to be a Calder\'{o}n-Zygmund standard
kernel.

Suppose that $K$ satisfies (\ref{1}) and (\ref{1*}). Then, Zhang and Wu \cite%
{Zhang} considered the family of operators $T:=\left\{ T_{\epsilon }\right\}
_{\epsilon >0}$ and a related the family of commutator operators $%
T_{b}:=\left\{ T_{\epsilon ,b}\right\} _{\epsilon >0}$ generated by $%
T_{\epsilon }$ and $b$ which are given by%
\begin{equation}
T_{\epsilon }f\left( x\right) =\dint\limits_{\left\vert x-y\right\vert
>\epsilon }K\left( x,y\right) f\left( y\right) dy  \label{3}
\end{equation}%
and%
\begin{equation}
T_{\epsilon ,b}f\left( x\right) =\dint\limits_{\left\vert x-y\right\vert
>\epsilon }\left( b\left( x\right) -b\left( y\right) \right) K\left(
x,y\right) f\left( y\right) dy.  \label{0}
\end{equation}%
In this sense, following \cite{Zhang}, the definition of the oscillation
operator of $T$ is given by

\begin{equation*}
\mathcal{O}\left( Tf\right) \left( x\right) :=\left(
\dsum\limits_{i=1}^{\infty }\sup_{t_{i+1}\leq \epsilon _{i+1}<\epsilon
_{i}\leq t_{i}}\left\vert T_{\epsilon _{i+1}}f\left( x\right) -T_{\epsilon
_{i}}f\left( x\right) \right\vert ^{2}\right) ^{\frac{1}{2}},
\end{equation*}%
where $\left\{ t_{i}\right\} $ is a decreasing fixed sequence of positive
numbers converging to $0$ and a related $\rho $-variation operator is
defined by%
\begin{equation*}
\mathcal{V}_{\rho }\left( Tf\right) \left( x\right) :=\sup_{\epsilon
_{i}\searrow 0}\left( \dsum\limits_{i=1}^{\infty }\left\vert T_{\epsilon
_{i+1}}f\left( x\right) -T_{\epsilon _{i}}f\left( x\right) \right\vert
^{\rho }\right) ^{\frac{1}{\rho }},\qquad \rho >2,
\end{equation*}%
where the supremum is taken over all sequences of real number $\left\{
\epsilon _{i}\right\} $ decreasing to $0$. We also take into account the
operator%
\begin{equation*}
\mathcal{O}^{\prime }\left( Tf\right) \left( x\right) :=\left(
\dsum\limits_{i=1}^{\infty }\sup_{t_{i+1}<\eta _{i}<t_{i}}\left\vert
T_{t_{i+1}}f\left( x\right) -T_{\eta _{i}}f\left( x\right) \right\vert
^{2}\right) ^{\frac{1}{2}}.
\end{equation*}%
On the other hand, it is obvious that 
\begin{equation*}
\mathcal{O}^{\prime }\left( Tf\right) \approx \mathcal{O}\left( Tf\right) .
\end{equation*}%
That is, 
\begin{equation*}
\mathcal{O}^{\prime }\left( Tf\right) \leq \mathcal{O}\left( Tf\right) \leq 2%
\mathcal{O}^{\prime }\left( Tf\right) .
\end{equation*}%
Recently, Campbell et al. in \cite{Campbell} proved the oscillation and
variation inequalities for the Hilbert transform in $L^{p}$($1<p<\infty $)
and then following \cite{Campbell}, we denote by $E$ the mixed norm Banach
space of two-variable function $h$ defined on $%
\mathbb{R}
\times 
\mathbb{N}
$ such that%
\begin{equation*}
\left\Vert h\right\Vert _{E}\equiv \left( \dsum\limits_{i}\left(
\sup_{s}\left\vert h\left( s,i\right) \right\vert \right) ^{2}\right)
^{1/2}<\infty .
\end{equation*}%
Given $T:=\left\{ T_{\epsilon }\right\} _{\epsilon >0}$ is a family
operators such that $\lim\limits_{\epsilon \rightarrow 0}T_{\epsilon
}f\left( x\right) =Tf\left( x\right) $ exists almost everywhere for certain
class of functions $f$, where $T_{\epsilon }$ defined as (\ref{3}). For a
fixed decreasing sequence $\left\{ t_{i}\right\} $ with $t_{i}\searrow 0$,
let $J_{i}=\left( t_{i+1},t_{i}\right] $ and define the $E$-valued operator $%
U\left( T\right) :f\rightarrow U\left( T\right) f$ given by%
\begin{equation*}
U\left( T\right) f\left( x\right) =\left\{ T_{t_{i+1}}f\left( x\right)
-T_{s}f\left( x\right) \right\} _{s\in J_{i},i\in 
\mathbb{N}
}=\left\{ \dint\limits_{\left\{ t_{i+1}<\left\vert x-y\right\vert <s\right\}
}K\left( x,y\right) f\left( y\right) dy\right\} _{s\in J_{i},i\in 
\mathbb{N}
}.
\end{equation*}%
Then 
\begin{eqnarray*}
\mathcal{O}^{\prime }\left( Tf\right) \left( x\right)  &=&\left\Vert U\left(
T\right) f\left( x\right) \right\Vert _{E}=\left\Vert \left\{
T_{t_{i+1}}f\left( x\right) -T_{s}f\left( x\right) \right\} _{s\in
J_{i},i\in 
\mathbb{N}
}\right\Vert _{E} \\
&=&\left\Vert \left\{ \dint\limits_{\left\{ t_{i+1}<\left\vert
x-y\right\vert <s\right\} }K\left( x,y\right) f\left( y\right) dy\right\}
_{s\in J_{i},i\in 
\mathbb{N}
}\right\Vert _{E}.
\end{eqnarray*}%
Let $\Phi =\left\{ \beta :\beta =\left\{ \epsilon _{i}\right\} ,\epsilon
_{i}\in 
\mathbb{R}
,\epsilon _{i}\searrow 0\right\} $. We denote by $F_{\rho }$ the mixed norm
space of two variable functions $g\left( i,\beta \right) $ such that 
\begin{equation*}
\left\Vert g\right\Vert _{F_{\rho }}\equiv \sup_{\beta }\left(
\dsum\limits_{i}\left\vert g\left( i,\beta \right) \right\vert ^{\rho
}\right) ^{1/\rho }.
\end{equation*}%
We also take into account the $F_{\rho }$-valued operator $V\left( T\right)
:f\rightarrow V\left( T\right) f$ such that 
\begin{equation*}
V\left( T\right) f\left( x\right) =\left\{ T_{\epsilon _{i+1}}f\left(
x\right) -T_{\epsilon _{i}}f\left( x\right) \right\} _{\beta =\left\{
\epsilon _{i}\right\} \in \Phi }.
\end{equation*}%
Thus,%
\begin{equation*}
V_{\rho }\left( T\right) f\left( x\right) =\left\Vert V\left( T\right)
f\left( x\right) \right\Vert _{F_{\rho }}.
\end{equation*}

Given $m$ is a positive integer, and $b$ is a function on ${\mathbb{R}}$.
Let $R_{m+1}\left( b;x,y\right) $ be the $m+1$-th order Taylor series
remainder of $b$ at $x$ about $y$, that is, 
\begin{equation*}
R_{m+1}\left( b;x,y\right) =b\left( x\right) -\dsum\limits_{\gamma \leq m}%
\frac{1}{\gamma !}b^{\left( \gamma \right) }\left( y\right) \left(
x-y\right) ^{\gamma }.
\end{equation*}%
In this paper, we consider the family of operators $T^{b}:=\left\{
T_{\epsilon }^{b}\right\} _{\epsilon >0}$ given by \cite{Hu}, where $%
T_{\epsilon }^{b}$ are the multilinear singular integral operators of $%
T_{\epsilon }$ as follows%
\begin{equation}
T_{\epsilon }^{b}f\left( x\right) =\dint\limits_{\left\vert x-y\right\vert
>\epsilon }\frac{R_{m+1}\left( b;x,y\right) }{\left\vert x-y\right\vert ^{m}}%
K\left( x,y\right) f\left( y\right) dy.  \label{4}
\end{equation}%
Thus, if $m=0$, then $T_{\epsilon }^{b}$ is just the commutator of $%
T_{\epsilon }$ and $b$, which is given by (\ref{0}). But, if $m>0$, then $%
T_{\epsilon }^{b}$ are non-trivial generation of the commutators.

The theory of multilinear analysis was received extensive studies in the
last 3 decades (see \cite{Cohen, Gurbuz} for example). Hu and Wang \cite{Hu}
proved that the weighted $\left( L^{p},L^{q}\right) $-boundedness of the
oscillation and variation operators for $T^{b}$ when the $m$-th derivative
of $b$ belongs to the homogenous Lipschitz space $\dot{\Lambda}_{\beta }$.
In this sense, we recall the definition of homogenous Lipschitz space $\dot{%
\Lambda}_{\beta }$ as follows:

\begin{definition}
$\left( \text{\textbf{Homogenous Lipschitz space}}\right) $ Let $0<\beta
\leq 1$. The homogeneous Lipschitz space $\dot{\Lambda}_{\beta }$ is defined
by%
\begin{equation*}
\dot{\Lambda}_{\beta }\left( {\mathbb{R}}\right) =\left\{ b:\left\Vert
b\right\Vert _{\dot{\Lambda}_{\beta }}=\sup_{x,h\in 
\mathbb{R}
,h\neq 0}\frac{\left\vert b\left( x+h\right) -b\left( x\right) \right\vert }{%
\left\vert h\right\vert ^{\beta }}<\infty \right\} .
\end{equation*}

Obviously, if $\beta >1$, then $\dot{\Lambda}_{\beta }\left( {\mathbb{R}}%
\right) $ only includes constant. So we restrict $0<\beta \leq 1$.
\end{definition}

Now, we recall the definitions of basic spaces such as Morrey, weighted
Lebesgue, weighted Morrey spaces and consider the relationship between these
spaces.

Besides the Lebesgue space $L^{q}\left( {\mathbb{R}}\right) $, the Morrey
space $M_{p}^{q}\left( {\mathbb{R}}\right) $ is another important function
space with definition as follows:

\begin{definition}
$\left( \text{\textbf{Morrey space}}\right) $\label{Definition1} For $1\leq
p\leq q<\infty $, the Morrey space $M_{p}^{q}\left( {\mathbb{R}}\right) $ is
the collection of all measurable functions $f$ whose Morrey space norm is%
\begin{equation*}
\left\Vert f\right\Vert _{M_{p}^{q}\left( {\mathbb{R}}\right) }=\sup 
_{\substack{ I\subset {\mathbb{R}}  \\ I:Interval}}\frac{1}{\left\vert
I\right\vert ^{\frac{1}{p}-\frac{1}{q}}}\left\Vert f\chi _{I}\right\Vert
_{L_{p}\left( {\mathbb{R}}\right) }<\infty .
\end{equation*}
\end{definition}

\begin{remark}
$\cdot $ If $p=q$, then%
\begin{equation*}
\Vert f\Vert _{M_{q}^{q}\left( {\mathbb{R}}\right) }=\Vert f\Vert
_{L^{q}\left( {\mathbb{R}}\right) }.
\end{equation*}

$\cdot $ if $q<p$, then $M_{p}^{q}\left( {\mathbb{R}}\right) $ is strictly
larger than $L^{q}\left( {\mathbb{R}}\right) $. For example, $%
f(x):=\left\vert x\right\vert ^{-\frac{1}{q}}\in $ $M_{p}^{q}\left( {\mathbb{%
R}}\right) $ but $f(x):=\left\vert x\right\vert ^{-\frac{1}{q}}\notin $ $%
L^{q}\left( {\mathbb{R}}\right) $.
\end{remark}

\bigskip On the other hand, for a given weight function $w$ and any interval 
$I$, we also denote the Lebesgue measure of $I$ by $\left\vert I\right\vert $
and set weighted measure 
\begin{equation*}
w\left( I\right) =\dint\limits_{I}w\left( x\right) dx.
\end{equation*}%
For $0<p<\infty $, the weighted Lebesgue space $L_{p}(w)\equiv L_{p}({{%
\mathbb{R}}},w)$ is defined by the norm 
\begin{equation*}
\Vert f\Vert _{L_{p}(w)}=\left( \dint\limits_{{{\mathbb{R}}}%
}|f(x)|^{p}w(x)dx\right) ^{\frac{1}{p}}<\infty .
\end{equation*}

A weight $w$ is said to belong to the Muckenhoupt class $A_{p}$ for $%
1<p<\infty $ such that%
\begin{align}
\lbrack w]_{A_{p}}& :=\sup\limits_{I}[w]_{A_{p}(I)}  \notag \\
& =\sup\limits_{I}\left( \frac{1}{|I|}\dint\limits_{I}w(x)dx\right) \left( 
\frac{1}{|I|}\dint\limits_{I}w(x)^{1-p^{\prime }}dx\right) ^{p-1}<\infty ,
\label{2}
\end{align}%
where $p^{\prime }=\frac{p}{p-1}$. The condition (\ref{2}) is called the $%
A_{p}$-condition, and the weights which satisfy it are called $A_{p}$%
-weights. The expression $[w]_{A_{p}}$ is also called characteristic
constant of $w$.

Here and after, $A_{p}$ denotes the Muckenhoupt classes (see \cite{Gurbuz,
KomShir}). The $A_{p}$ class of weights characterizes the $L_{p}(w)$
boundedness of the maximal function as Muckenhoupt \cite{Muckenhoupt}
established in the 70s. Subsequent works of Muckenhoupt \cite{Muckenhoupt}
himself Muckenhoupt and Wheeden \cite{Muckenhoupt1, Muckenhoupt2}, Coifman
and Fefferman \cite{Coifman} were devoted to explore the connection of the $%
A_{p}$ class with weighted estimates for singular integrals. However, it was
not until the 2000s that the quantitative dependence on the so called $A_{p}$
constant, namely $[w]_{A_{p}}$, became a trending topic.

When $p=1$, $w\in $ $A_{1}$ if there exists $C>1$ such that for almost every 
$x$, 
\begin{equation}
Mw(x)dx\leq Cw\left( x\right)  \label{5}
\end{equation}%
and the infimum of $C$ satisfying the inequality (\ref{5}) is denoted by $%
[w]_{A_{1}}$, where $M$ is the classical Hardy-Littlewood maximal operator.

When $p=\infty $, we define $A_{\infty }\left( {{\mathbb{R}}}\right)
=\dbigcup\limits_{1\leq p<\infty }A_{p}\left( {{\mathbb{R}}}\right) $. That
is, the $A_{\infty }$ constant is given by%
\begin{eqnarray*}
\lbrack w]_{A_{\infty }} &:&=\sup\limits_{I}[w]_{A_{\infty }(I)} \\
&=&\sup\limits_{I}\dint\limits_{I}M\left( \chi _{I}w\right) \left( x\right)
dx,
\end{eqnarray*}%
where we utilize the notation $M\left( \chi _{I}w\right) $ to denote the
Hardy-Littlewood maximal function of a function $\chi _{I}w$ by%
\begin{equation*}
M\left( \chi _{I}w\right) (x):=\sup\limits_{I}\frac{1}{|I|}%
\int\limits_{I}|\chi _{I}w(x)|dx.
\end{equation*}

\bigskip A weight function $w$ belongs to $A_{p,q}$ (Muckenhoupt-Wheeden
class) for $1<p<q<\infty $ if 
\begin{align}
\lbrack w]_{A_{p,q}}& :=\sup\limits_{I}[w]_{A_{p,q}(I)}  \notag \\
& =\sup\limits_{I}\left( \frac{1}{|I|}\dint\limits_{I}w(x)^{q}dx\right) ^{%
\frac{1}{q}}\left( \frac{1}{|I|}\dint\limits_{I}w(x)^{-p^{\prime }}dx\right)
^{\frac{1}{p^{\prime }}}<\infty .  \label{13}
\end{align}%
From the definition of $A_{p,q}$, we know that $w\left( x\right) \in
A_{p,q}\left( {{\mathbb{R}}}\right) $ implies $w(x)^{q}\in A_{q}\left( {{%
\mathbb{R}}}\right) $ and $w(x)^{p}\in A_{p}\left( {{\mathbb{R}}}\right) $.

Now, we begin with some Lemmas. These Lemmas are very necessary for the
proof of the main result.

\begin{lemma}
\label{Lemma1}\cite{GarRub} If $w\in A_{p}$, $p\geq 1$, then there exists a
constant $C>0$ such that%
\begin{equation*}
w\left( 2I\right) \leq Cw\left( I\right) .
\end{equation*}%
for any interval $I$.

More precisely, for all $\lambda >1$ we have%
\begin{equation*}
w\left( \lambda I\right) \leq C\lambda ^{p}w\left( I\right) ,
\end{equation*}%
where $C$ is a constant independent of $I$ or $\lambda $ and $w\left(
I\right) =\dint\limits_{I}w\left( x\right) dx$.
\end{lemma}

\begin{lemma}
\label{Lemma2}\cite{Cohen} Let $b$ be a function on $%
\mathbb{R}
$ and $b^{\left( m\right) }\in L_{u}\left( 
\mathbb{R}
\right) $ with $m\in 
\mathbb{N}
$ for any $u>1$. Then%
\begin{equation*}
\left\vert R_{m}\left( b;x,y\right) \right\vert \leq C\left\vert
x-y\right\vert ^{m}\left( \frac{1}{\left\vert I\left( x,y\right) \right\vert 
}\dint\limits_{I\left( x,y\right) }\left\vert b^{\left( m\right) }\left(
z\right) \right\vert ^{u}dz\right) ^{\frac{1}{u}},C>0,
\end{equation*}%
where $I\left( x,y\right) $ is the interval $\left( x-5\left\vert
x-y\right\vert ,x+5\left\vert x-y\right\vert \right) $.
\end{lemma}

\begin{lemma}
\label{Lemma3}\cite{Hu} Let $K\left( x,y\right) $ satisfies (\ref{1}) and (%
\ref{1*}), $\rho >2$, and $T:=\left\{ T_{\epsilon }\right\} _{\epsilon >0}$
and $T^{b}:=\left\{ T_{\epsilon }^{b}\right\} _{\epsilon >0}$ be given by (%
\ref{3}) and (\ref{4}), respectively. If $\mathcal{O}\left( T\right) $ and $%
\mathcal{V}_{\rho }\left( T\right) $ are bounded on $L_{p_{0}}\left( 
\mathbb{R}
,dx\right) $ for some $1<p_{0}<\infty $, and $b^{\left( m\right) }\in \dot{%
\Lambda}_{\beta }$ with $m\in 
\mathbb{N}
$ for $0<\beta <1$, then 
\begin{equation}
\left\Vert \mathcal{O}^{\prime }\left( T^{b}\right) \right\Vert
_{L_{q}\left( w^{q}\right) }\leq \left\Vert \mathcal{O}\left( T^{b}\right)
\right\Vert _{L_{q}\left( w^{q}\right) }\leq C\left\Vert b\right\Vert _{\dot{%
\Lambda}_{\beta }}\Vert f\Vert _{L_{p}(w^{p})},C>0,  \label{6}
\end{equation}%
and 
\begin{equation*}
\left\Vert \mathcal{V}_{\rho }\left( T^{b}\right) \right\Vert _{L_{q}\left(
w^{q}\right) }\leq C\left\Vert b\right\Vert _{\dot{\Lambda}_{\beta }}\Vert
f\Vert _{L_{p}(w^{p})},C>0,
\end{equation*}
for any $1<p<\frac{1}{\beta }$ with $\frac{1}{q}=\frac{1}{p}-\beta $ and $%
w\in A_{p,q}$.
\end{lemma}

Next, in 2009, the weighted Morrey space $L_{p,\kappa }(w)$ was defined by
Komori and Shirai \cite{KomShir} as follows:

\begin{definition}
$\left( \text{\textbf{Weighted Morrey space}}\right) $ Let $1\leq p<\infty $%
, $0<\kappa <1$ and $w$ be a weight function. Then the weighted Morrey space 
$L_{p,\kappa }(w)\equiv L_{p,\kappa }({\mathbb{R}},w)$ is defined by 
\begin{equation*}
L_{p,\kappa }(w)\equiv L_{p,\kappa }({\mathbb{R}},w)=\left\{ f\in
L_{p,w}^{loc}\left( {\mathbb{R}}\right) :\Vert f\Vert _{L_{p,\kappa
}(w)}=\sup\limits_{I}\,w(I)^{-\frac{\kappa }{p}}\,\Vert f\Vert
_{L_{p,w}(I)}<\infty \right\} .
\end{equation*}
\end{definition}

\begin{remark}
$\cdot $ If $\kappa =0,$ then%
\begin{equation*}
\Vert f\Vert _{L_{p,0}(w)}=\Vert f\Vert _{L_{p}(w)}.
\end{equation*}

$\cdot $ When $w\equiv 1$ and $\kappa =1-\frac{p}{q}$ with $1<p\leq q<\infty 
$, then%
\begin{equation*}
\Vert f\Vert _{L_{p,1-\frac{p}{q}}(1)}=\Vert f\Vert _{M_{p}^{q}\left( {%
\mathbb{R}}\right) }.
\end{equation*}
\end{remark}

Finally, we recall the definition of the weighted Morrey space with two
weights as follows:

\begin{definition}
$\left( \text{\textbf{Weighted Morrey space with two weights}}\right) $ Let $%
1\leq p<\infty $ and $0<\kappa <1$. Then for two weights $u$ and $v$, the
weighted Morrey space $L_{p,\kappa }(u,v)\equiv L_{p,\kappa }({\mathbb{R}}%
,u,v)$ is defined by%
\begin{equation*}
L_{p,\kappa }(u,v)\equiv L_{p,\kappa }({\mathbb{R}},u,v)=\left\{ f\in
L_{p,u}^{loc}\left( {\mathbb{R}}\right) :\Vert f\Vert _{L_{p,\kappa
}(w)}=\sup\limits_{I}\,v(I)^{-\frac{\kappa }{p}}\,\Vert f\Vert
_{L_{p,u}(I)}<\infty \right\} .
\end{equation*}%
It is obvious that%
\begin{equation*}
L_{p,\kappa }(w,w)\equiv L_{p,\kappa }(w).
\end{equation*}
\end{definition}

In 2016, Zhang and Wu \cite{Zhang} gave the boundedness of the oscillation
and variation operators for Calder\'{o}n-Zygmund singular integrals and the
corresponding commutators on the weighted Morrey spaces. In 2017, Hu and
Wang \cite{Hu} established the weighted $\left( L^{p},L^{q}\right) $%
-inequalities of the variation and oscillation operators for the multilinear
Calder\'{o}n-Zygmund singular integral with a Lipschitz function in $%
\mathbb{R}
$. Inspired of these results \cite{Hu, Zhang}, we investigate the
boundedness of the oscillation and variation operators for the family of the
multilinear singular integral defined by (\ref{4}) on weighted Morrey spaces
when the $m$-th derivative of $b$ belongs to the homogenous Lipschitz space $%
\dot{\Lambda}_{\beta }$ in this work.

Throughout this paper, $C$ always means a positive constant independent of
the main parameters involved, and may change from one occurrence to another.
We also use the notation $F\lesssim G$ to mean $F\leq CG$ for an appropriate
constant $C>0$, and $F\approx G$ to mean $F\lesssim G$ and $G\lesssim F$.

\section{Main result}

We now formulate our main result as follows.

\begin{theorem}
\label{Theorem1}Let $K\left( x,y\right) $ satisfies (\ref{1}) and (\ref{1*}%
), $\rho >2$, and $T:=\left\{ T_{\epsilon }\right\} _{\epsilon >0}$ and $%
T^{b}:=\left\{ T_{\epsilon }^{b}\right\} _{\epsilon >0}$ be given by (\ref{3}%
) and (\ref{4}), respectively. If $\mathcal{O}\left( T\right) $ and $%
\mathcal{V}_{\rho }\left( T\right) $ are bounded on $L_{p_{0}}\left( 
\mathbb{R}
,dx\right) $ for some $1<p_{0}<\infty $, and $b^{\left( m\right) }\in \dot{%
\Lambda}_{\beta }$ with $m\in 
\mathbb{N}
$ for $0<\beta <1$, then $\mathcal{O}\left( T^{b}\right) $ and $\mathcal{V}%
_{\rho }\left( T^{b}\right) $ are bounded from $L_{p,\kappa }(w^{p},w^{q})$
to $L_{p,\frac{\kappa q}{p}}(w^{q})$ for any $1<p<\frac{1}{\beta }$, $\frac{1%
}{q}=\frac{1}{p}-\beta $, $0<\kappa <\frac{p}{q}$ and $w\in A_{p,q}$.
\end{theorem}

\begin{corollary}
\cite{Zhang} Let $K\left( x,y\right) $ satisfies (\ref{1}) and (\ref{1*}), $%
\rho >2$, and $T:=\left\{ T_{\epsilon }\right\} _{\epsilon >0}$ and $%
T_{b}:=\left\{ T_{\epsilon ,b}\right\} _{\epsilon >0}$ be given by (\ref{3})
and (\ref{0}), respectively. If $\mathcal{O}\left( T\right) $ and $\mathcal{V%
}_{\rho }\left( T\right) $ are bounded on $L_{p_{0}}\left( 
\mathbb{R}
,dx\right) $ for some $1<p_{0}<\infty $, and $b\in \dot{\Lambda}_{\beta }$
for $0<\beta <1$, then $\mathcal{O}\left( T_{b}\right) $ and $\mathcal{V}%
_{\rho }\left( T_{b}\right) $ are bounded from $L_{p,\kappa }(w^{p},w^{q})$
to $L_{p,\frac{\kappa q}{p}}(w^{q})$ for any $1<p<\frac{1}{\beta }$, $\frac{1%
}{q}=\frac{1}{p}-\beta $, $0<\kappa <\frac{p}{q}$ and $w\in A_{p,q}$.
\end{corollary}

\subsection{The Proof of Theorem \protect\ref{Theorem1}}

\begin{proof}
We consider the proof related to $\mathcal{O}\left( T^{b}\right) $ firstly.
Fix an interval $I=\left( x_{0}-l,x_{0}+l\right) $, and we write as $%
f=f_{1}+f_{2}$, where $f_{1}=f\chi _{2I}$, $\chi _{2I}$ denotes the
characteristic function of $2I$. Thus, it is sufficient to show that the
conclusion%
\begin{eqnarray*}
\left\Vert \mathcal{O}^{\prime }\left( T^{b}f\right) \left( x\right)
\right\Vert _{L_{p,\frac{\kappa q}{p}}(w^{q})} &\leq &\left\Vert \mathcal{O}%
^{\prime }\left( T^{b}f_{1}\right) \left( x\right) \right\Vert _{L_{p,\frac{%
\kappa q}{p}}(w^{q})}+\left\Vert \mathcal{O}^{\prime }\left(
T^{b}f_{1}\right) \left( x\right) \right\Vert _{L_{p,\frac{\kappa q}{p}%
}(w^{q})} \\
&\lesssim &\left\Vert b\right\Vert _{\dot{\Lambda}_{\beta }}\left\Vert
f\right\Vert _{L_{p,\kappa }(w^{p},w^{q})}
\end{eqnarray*}%
holds for every interval $I\subset {\mathbb{R}}$. Then%
\begin{eqnarray*}
&&\left( \dint\limits_{I}\left\vert \mathcal{O}^{\prime }\left(
T^{b}f\right) \left( x\right) \right\vert ^{q}w^{q}\left( x\right) dx\right)
^{\frac{1}{q}} \\
&\leq &\left( \dint\limits_{I}\left\vert \mathcal{O}^{\prime }\left(
T^{b}f_{1}\right) \left( x\right) \right\vert ^{q}w^{q}\left( x\right)
dx\right) ^{\frac{1}{q}}+\left( \dint\limits_{I}\left\vert \mathcal{O}%
^{\prime }\left( T^{b}f_{2}\right) \left( x\right) \right\vert
^{q}w^{q}\left( x\right) dx\right) ^{\frac{1}{q}} \\
&=&:F_{1}+F_{2}.
\end{eqnarray*}

First, we use (\ref{6}) to estimate $F_{1}$, and we obtain%
\begin{eqnarray*}
F_{1} &=&\left( \dint\limits_{I}\left\vert \mathcal{O}^{\prime }\left(
T^{b}f_{1}\right) \left( x\right) \right\vert ^{q}w^{q}\left( x\right)
dx\right) ^{\frac{1}{q}}\lesssim \left\Vert b\right\Vert _{\dot{\Lambda}%
_{\beta }}\Vert f_{1}\Vert _{L_{p}(w^{p})} \\
&=&\left\Vert b\right\Vert _{\dot{\Lambda}_{\beta }}\left( \frac{1}{%
w^{q}\left( 2I\right) ^{\kappa }}\dint\limits_{2I}\left\vert f\left(
x\right) \right\vert ^{p}w^{p}\left( x\right) dx\right) ^{\frac{1}{p}%
}w^{q}\left( 2I\right) ^{\frac{\kappa }{p}} \\
&\lesssim &\left\Vert b\right\Vert _{\dot{\Lambda}_{\beta }}\left\Vert
f\right\Vert _{L_{p,\kappa }(w^{p},w^{q})}^{p}w^{q}\left( I\right) ^{\frac{%
\kappa }{p}}.
\end{eqnarray*}%
Thus, 
\begin{equation}
\left\Vert \mathcal{O}^{\prime }\left( T^{b}f_{1}\right) \left( x\right)
\right\Vert _{L_{p,\frac{\kappa q}{p}}(w^{q})}\lesssim \left\Vert
b\right\Vert _{\dot{\Lambda}_{\beta }}\left\Vert f\right\Vert _{L_{p,\kappa
}(w^{p},w^{q})}.  \label{20}
\end{equation}%
Second, for $x\in I,k=1,2,\ldots ,m\in 
\mathbb{N}
$, let $A_{k}=\left\{ y:2^{k}l\leq \left\vert y-x\right\vert
<2^{k+1}l\right\} $, $B_{k}=\left\{ y:\left\vert y-x\right\vert
<2^{k+1}l\right\} $, and 
\begin{equation*}
b_{k}\left( z\right) =b\left( z\right) -\frac{1}{m!}\left( b^{\left(
m\right) }\right) _{B_{k}}z^{m}.
\end{equation*}%
By \cite{Cohen}, for any $y\in A_{k}$, it is obvious that 
\begin{equation*}
R_{m+1}\left( b;x,y\right) =R_{m+1}\left( b_{k};x,y\right) .
\end{equation*}%
Moreover, since $b\in \dot{\Lambda}_{\beta }$, then, for $y\in A_{k}$, we get%
\begin{eqnarray}
\left\vert b^{\left( m\right) }\left( y\right) -\left( b^{\left( m\right)
}\right) _{B_{k}}\right\vert &\leq &\frac{1}{\left\vert B_{k}\right\vert }%
\dint\limits_{B_{k}}\left\vert b^{\left( m\right) }\left( y\right)
-b^{\left( m\right) }\left( z\right) \right\vert dz  \notag \\
&\lesssim &\left\Vert b^{\left( m\right) }\right\Vert _{\dot{\Lambda}_{\beta
}}\left( 2^{k}l\right) ^{\beta }.  \label{7}
\end{eqnarray}%
Hence, by Lemma \ref{Lemma2} and (\ref{7})%
\begin{eqnarray*}
R_{m}\left( b_{k};x,y\right) &\lesssim &\left\vert x-y\right\vert ^{m}\left( 
\frac{1}{\left\vert I\left( x,y\right) \right\vert }\dint\limits_{I\left(
x,y\right) }\left\vert b^{\left( m\right) }\left( z\right) \right\vert
^{u}dz\right) ^{\frac{1}{u}} \\
&\lesssim &\left\vert x-y\right\vert ^{m}\left\Vert b^{\left( m\right)
}\right\Vert _{\dot{\Lambda}_{\beta }}\left( 2^{k}l\right) ^{\beta }.
\end{eqnarray*}%
Also, following \cite{Zhang}, we have%
\begin{equation*}
\left\Vert \left\{ \chi _{\left\{ t_{i}+1<\left\vert x-y\right\vert
<u\right\} }\right\} _{u\in J_{i},i\in 
\mathbb{N}
}\right\Vert _{A}\leq 1.
\end{equation*}%
Thus, the estimate of $F_{2}$ can be obtained as follows:%
\begin{eqnarray*}
\left\vert \mathcal{O}^{\prime }\left( T^{b}f_{2}\right) \left( x\right)
\right\vert &=&\left\Vert U\left( T^{b}f_{2}\right) \left( x\right)
\right\Vert \\
&=&\left\Vert \left\{ \dint\limits_{\left\{ t_{i}+1<\left\vert
x-y\right\vert <u\right\} }\frac{R_{m+1}\left( b;x,y\right) }{\left\vert
x-y\right\vert ^{m}}K\left( x,y\right) f_{2}\left( y\right) dy\right\}
\right\Vert _{A} \\
&\leq &\dint\limits_{%
\mathbb{R}
}\left\Vert \left\{ \chi _{\left\{ t_{i}+1<\left\vert x-y\right\vert
<u\right\} }\right\} _{u\in J_{i},i\in 
\mathbb{N}
}\right\Vert _{A}\left\vert \frac{R_{m+1}\left( b;x,y\right) }{\left\vert
x-y\right\vert ^{m}}K\left( x,y\right) f_{2}\left( y\right) \right\vert dy \\
&\leq &\dint\limits_{%
\mathbb{R}
}\left\vert \frac{R_{m+1}\left( b;x,y\right) }{\left\vert x-y\right\vert ^{m}%
}K\left( x,y\right) f_{2}\left( y\right) \right\vert dy \\
&\lesssim &\dint\limits_{\left\vert x-y\right\vert >2l}\left\vert \frac{%
R_{m+1}\left( b;x,y\right) }{\left\vert x-y\right\vert ^{m}}K\left(
x,y\right) f\left( y\right) \right\vert dy \\
&\lesssim &\dsum\limits_{k=1}^{\infty }\frac{1}{2^{k}l}\dint\limits_{A_{k}}%
\left( \left\Vert b^{\left( m\right) }\right\Vert _{\dot{\Lambda}_{\beta
}}\left( 2^{k}l\right) ^{\beta }+\left\vert b^{\left( m\right) }\left(
y\right) -\left( b^{\left( m\right) }\right) _{B_{k}}\right\vert \right)
\left\vert f\left( y\right) \right\vert dy \\
&\lesssim &\left\Vert b^{\left( m\right) }\right\Vert _{\dot{\Lambda}_{\beta
}}\dsum\limits_{k=1}^{\infty }\frac{1}{\left( 2^{k}l\right) ^{1-\beta }}%
\dint\limits_{A_{k}}\left\vert f\left( y\right) \right\vert
dy+\dsum\limits_{k=1}^{\infty }\frac{1}{2^{k}l}\dint\limits_{A_{k}}\left%
\vert b^{\left( m\right) }\left( y\right) -\left( b^{\left( m\right)
}\right) _{B_{k}}\right\vert \left\vert f\left( y\right) \right\vert dy \\
&=&G_{1}+G_{2}.
\end{eqnarray*}%
For $G_{1}$, since%
\begin{equation*}
\left( \dint\limits_{A_{k}}w\left( y\right) ^{-p^{\prime }}dy\right) ^{\frac{%
1}{p^{\prime }}}\lesssim w^{q}\left( B_{k}\right) ^{-\frac{1}{q}}\left\vert
B_{k}\right\vert ^{\frac{1}{p^{\prime }}+\frac{1}{q}}
\end{equation*}%
with $1<p<\frac{1}{\beta }$, $\frac{1}{q}=\frac{1}{p}-\beta $ and using H%
\"{o}lder's inequality, we have%
\begin{eqnarray}
&&\dsum\limits_{k=1}^{\infty }\frac{1}{\left( 2^{k}l\right) ^{1-\beta }}%
\dint\limits_{A_{k}}\left\vert f\left( y\right) \right\vert dy  \notag \\
&\lesssim &\dsum\limits_{k=1}^{\infty }\frac{1}{\left( 2^{k}l\right)
^{1-\beta }}\left( \dint\limits_{A_{k}}\left\vert f\left( y\right)
\right\vert ^{p}w^{p}\left( y\right) dy\right) ^{\frac{1}{p}}\left(
\dint\limits_{A_{k}}w\left( y\right) ^{-p^{\prime }}dy\right) ^{\frac{1}{%
p^{\prime }}}  \notag \\
&\lesssim &\left\Vert f\right\Vert _{L_{p,\kappa
}(w^{p},w^{q})}\dsum\limits_{k=1}^{\infty }\frac{\left( 2^{k}l\right) ^{%
\frac{1}{p^{\prime }}+\frac{1}{q}}}{\left( 2^{k}l\right) ^{1-\beta }}%
w^{q}\left( B_{k}\right) ^{\frac{\kappa }{p}-\frac{1}{q}}  \notag \\
&\lesssim &\left\Vert f\right\Vert _{L_{p,\kappa
}(w^{p},w^{q})}\dsum\limits_{k=1}^{\infty }w^{q}\left( B_{k}\right) ^{\frac{%
\kappa }{p}-\frac{1}{q}}.  \label{10}
\end{eqnarray}%
Since $w\in A_{p,q}$, then we have $w^{q}\in A_{\infty }$. Thus, Lemma \ref%
{Lemma1} implies $w^{q}\left( B_{k}\right) \leq \left( C\right)
^{k}w^{q}\left( I\right) ,C>1,$ i.e.,%
\begin{equation}
\dsum\limits_{k=1}^{\infty }w^{q}\left( B_{k}\right) ^{\frac{\kappa }{p}-%
\frac{1}{q}}\lesssim w^{q}\left( I\right) ^{\frac{\kappa }{p}-\frac{1}{q}%
}\dsum\limits_{k=1}^{\infty }C^{\frac{\kappa }{p}-\frac{1}{q}}\lesssim
w^{q}\left( I\right) ^{\frac{\kappa }{p}-\frac{1}{q}}  \label{11}
\end{equation}%
with $\frac{\kappa }{p}-\frac{1}{q}<0$. This implies%
\begin{equation}
G_{1}\lesssim \left\Vert b^{\left( m\right) }\right\Vert _{\dot{\Lambda}%
_{\beta }}\left\Vert f\right\Vert _{L_{p,\kappa }(w^{p},w^{q})}w^{q}\left(
I\right) ^{\frac{\kappa }{p}-\frac{1}{q}}.  \label{12}
\end{equation}%
Let $y\in A_{k}$. For $G_{2}$, by (\ref{7}), (\ref{10}) and (\ref{11}) we get%
\begin{eqnarray}
G_{2} &\lesssim &\left\Vert b^{\left( m\right) }\right\Vert _{\dot{\Lambda}%
_{\beta }}\dsum\limits_{k=1}^{\infty }\frac{1}{\left( 2^{k+1}l\right)
^{1-\beta }}\dint\limits_{A_{k}}\left\vert f\left( y\right) \right\vert dy 
\notag \\
&\lesssim &\left\Vert b^{\left( m\right) }\right\Vert _{\dot{\Lambda}_{\beta
}}\left\Vert f\right\Vert _{L_{p,\kappa }(w^{p},w^{q})}w^{q}\left( I\right)
^{\frac{\kappa }{p}-\frac{1}{q}}.  \label{12*}
\end{eqnarray}%
Thus, by (\ref{12}) and (\ref{12*}), we obtain 
\begin{eqnarray*}
F_{2} &=&\left( \dint\limits_{I}\left\vert \mathcal{O}^{\prime }\left(
T^{b}f_{2}\right) \left( x\right) \right\vert ^{q}w^{q}\left( x\right)
dx\right) ^{\frac{1}{q}} \\
&\lesssim &\left\Vert b^{\left( m\right) }\right\Vert _{\dot{\Lambda}_{\beta
}}\left\Vert f\right\Vert _{L_{p,\kappa }(w^{p},w^{q})}w^{q}\left( I\right)
^{\frac{\kappa }{p}-\frac{1}{q}}w^{q}\left( I\right) ^{\frac{1}{q}} \\
&=&\left\Vert b^{\left( m\right) }\right\Vert _{\dot{\Lambda}_{\beta
}}\left\Vert f\right\Vert _{L_{p,\kappa }(w^{p},w^{q})}w^{q}\left( I\right)
^{\frac{\kappa }{p}}.
\end{eqnarray*}%
Thus, 
\begin{equation}
\left\Vert \mathcal{O}^{\prime }\left( T^{b}f_{2}\right) \left( x\right)
\right\Vert _{L_{p,\frac{\kappa q}{p}}(w^{q})}\lesssim \left\Vert
b\right\Vert _{\dot{\Lambda}_{\beta }}\left\Vert f\right\Vert _{L_{p,\kappa
}(w^{p},w^{q})}.  \label{21}
\end{equation}%
As a result, by (\ref{20}) and (\ref{21}), we get%
\begin{equation*}
\left\Vert \mathcal{O}^{\prime }\left( T^{b}f\right) \left( x\right)
\right\Vert _{L_{p,\frac{\kappa q}{p}}(w^{q})}\lesssim \left\Vert
b\right\Vert _{\dot{\Lambda}_{\beta }}\left\Vert f\right\Vert _{L_{p,\kappa
}(w^{p},w^{q})}.
\end{equation*}%
Similarly, $\mathcal{V}_{\rho }\left( T^{b}\right) $ has the same estimate
as above (here we omit the details), thus the inequality 
\begin{equation*}
\left\Vert \mathcal{V}_{\rho }\left( T^{b}f\right) \left( x\right)
\right\Vert _{L_{p,\frac{\kappa q}{p}}(w^{q})}\lesssim \left\Vert
b\right\Vert _{\dot{\Lambda}_{\beta }}\left\Vert f\right\Vert _{L_{p,\kappa
}(w^{p},w^{q})}
\end{equation*}%
is valid.

Therefore, Theorem \ref{Theorem1} is completely proved.
\end{proof}


\begin{thebibliography}{99}
\bibitem{Campbell} J.T. Campbell, R.L. Jones, K. Reinhold and M. Wierdl,
Oscillation and variation for the Hilbert transform, Duke Math. J., \textbf{%
105} (2000), 59-83.

\bibitem{Cohen} J. Cohen and J. Gosselin, A $BMO$\ estimate for multilinear
singular integrals\textit{. }Illinois J. Math.\ \textbf{30}(1986), 445-464.

\bibitem{Coifman} R.R. Coifman and C. Fefferman, Weighted norm inequalities
for maximal functions and singular integrals, Studia Mathematica 51 (3)
(1974), 241-250.

\bibitem{GarRub} J. Garcia-Cuerva, J.L. Rubio de Francia, \textit{Weighted
norm inequalities and related topics. }North-Holland Math.,\textit{\ }%
\textbf{116}, Amsterdam, 1985.

\bibitem{Gurbuz} F. G\"{u}rb\"{u}z, On the behaviors of rough multilinear
fractional integral{\ }and multi-sublinear fractional maximal{\ operators
both on product }$L^{p}$ and {\ weighted }$L^{p}$ spaces, arXiv:1603.03466v2
[math.CA] 16 Jan 2018.

\bibitem{Hu} Y. Hu and Y.S. Wang, Oscillation and variation inequalities for
the multilinear singular integrals related to Lipschitz functions J.
Inequal. Appl. (2017) \textbf{2017}:292, 1-14.

\bibitem{KomShir} Y. Komori, S. Shirai, Weighted Morrey spaces and a
singular integral operator, Math. Nachr., \textbf{282}(2) (2009), 219-231.

\bibitem{Morrey} C.B. Morrey, On the solutions of quasi-linear elliptic
partial differential equations, Trans. Amer. Math. Soc., \textbf{43} (1938),
126-166.

\bibitem{Muckenhoupt1} B. Muckenhoupt and R. L. Wheeden, Weighted norm
inequalities for singular and fractional integrals\textit{, }Trans. Amer.
Math. Soc. \textbf{161}(1971), 249-258.

\bibitem{Muckenhoupt} B. Muckenhoupt, Weighted norm inequalities for the
Hardy maximal function, Trans. Amer. Math. Soc., \textbf{165} (1972),
207-226.

\bibitem{Muckenhoupt2} B. Muckenhoupt and R. L. Wheeden, Weighted norm
inequalities for fractional integrals, Trans. Amer. Math. Soc., \textbf{192}
(1974), 261-274.

\bibitem{Zhang} J. Zhang and H.X. Wu, Oscillation and variation inequalities
for singular integrals and commutators on weighted Morrey spaces, Front.
Math. China, \textbf{11}(2) (2016), 423-447.
\end{thebibliography}
\end{document}